\documentclass[12pt,a4paper,twoside]{article}
\usepackage{latexsym}

\usepackage{amssymb}
\usepackage{amsmath}
\usepackage{a4wide}
\title{New Classes of Almost Bent and Almost Perfect Nonlinear Polynomials\footnote{Part of this paper was presented at WCC 2005 \cite{BCP}.}}
\author{Lilya Budaghyan\footnote{Institute of Algebra and Geometry,
    Otto-von-Guericke  University  Magdeburg, D-39016   Magdeburg, GERMANY;
    e-mail: lilya.budaghyan@student.uni-magdeburg.de.},
\and{Claude Carlet\footnote{INRIA, Projet CODES,
BP 105 - 78153, Le Chesnay Cedex, FRANCE; e-mail: claude.carlet@inria.fr;
also member of the University of Paris 8.},}
\and{ Alexander Pott\footnote{Institute of Algebra and Geometry, 
Otto-von-Guericke  University  Magdeburg, D-39016   Magdeburg, GERMANY; e-mail: 
alexander.pott@mathematik.uni-magdeburg.de.}}}

\begin{document}
\date{}
\maketitle

\begin{abstract}  We construct infinite classes of almost bent and almost 
perfect nonlinear polynomials, which are affinely inequivalent to any sum of a 
power function and an affine function.
\end{abstract}
\section{Introduction}
Vectorial Boolean functions are used in cryptography, more precisely in block 
ciphers. An important condition on these functions is a high resistance to the 
differential  and linear cryptanalyses, which are the main attacks on block 
ciphers.
The functions with the smallest possible differential uniformity (see 
\cite{nyb}) oppose an optimum resistance to the differential attack  
\cite{B-Sh}.  They are called  almost perfect nonlinear (APN). The functions 
achieving the maximal possible nonlinearity (see \cite{CV,nyb1}) possess the 
best resistance to the linear attack \cite{M} and they are called almost bent 
(AB) or maximum nonlinear.

Up to now, in the study of APN and AB functions the main attention has been 
payed to power mappings and  all known constructions of APN and AB functions 
happen to be equivalent to power functions (see for example 
\cite{CCD3,CCD2,CCD,D1,D2,D3}). Recall that functions $F$ and $F^\prime$ are 
called equivalent if either $F^\prime=A_1\circ F\circ A_2 +A$ or 
$F^\prime=A_1\circ F^{-1}\circ A_2 +A$ (in case $F$ is a permutation) for some 
affine functions $A$, $A_1$ and $A_2$, where $A_1$ and $A_2$ are permutations. 
We shall say that two functions $F$ and $F^\prime=A_1\circ F\circ A_2 +A$ (with 
$A,A_{1},A_{2}$ affine, $A_{1},A_{2}$ being permutations) are  extended affine 
equivalent (EA-equivalent).
In this paper we give the first theoretical constructions of APN and AB 
polynomials which are inequivalent to power mappings. In these constructions, we 
apply the transformation of functions given by Carlet, Charpin and Zinoviev (see 
Proposition 3 in \cite{CC-carchzin}) to the Gold APN and AB  mappings 
\cite{Gold,nyb}. This transformation leads to an equivalence relation of 
functions that we call the Carlet-Charpin-Zinoviev equivalence 
(CCZ-equivalence). CCZ-equivalence corresponds to the affine equivalence of the 
graphs of functions, i.e. functions $F$ and $F^\prime$ are CCZ-equivalent if and 
only if, for some affine permutation, the image of the graph of $F$ is the graph 
of the function $F^\prime$. CCZ-equivalence preserves the nonlinearity, the 
differential uniformity of functions and the resistance of a function to the 
algebraic cryptanalysis \cite{CourtPier}.
It could be expected that CCZ-equivalence coincides in practice with the 
equivalence mentioned above. The present paper aims at proving and illustrating 
that  CCZ-equivalence is more general. 

The next section contains all the necessary definitions related to vectorial 
Boolean functions, including EA-equivalence, APN and AB properties. 

In Section 3, we give the definition of CCZ-equivalence, we describe its main 
properties and we show its connections with EA-equivalence. 
 
We give some results related to a classification of functions CCZ-equivalent to 
the Gold mappings in Section 4. 

Theorems 1, 2 and 4 in Section 5 present constructions of AB and APN polynomials 
which are EA-inequivalent to power functions and Theorem 3 presents an APN 
function EA-inequivalent to all known APN mappings. 

This paper is an extended version of an abstract published in \cite{BCP}. In particular, Theorem 4, the proof of Theorem 3 and some properties of CCZ-equivalence are not included in the abstract.

\section{Almost perfect nonlinear and almost bent functions}
Let $\mathbb{F}_{2}^m$ be the $m$-dimensional vector space over the field 
$\mathbb{F}_{2}$. Any function $F$ from $\mathbb{F}_{2}^m$ to itself can be 
uniquely represented as a polynomial on $m$ variables with coefficients in 
$\mathbb{F}_{2}^m$, whose degree with respect to each coordinate is at most 1:
\begin{displaymath}
F(x_1,...,x_m)=\sum_{u\in \mathbb{F}_{2}^m}c(u)\big(\prod_{i=1}^mx_i^{u_i}\big), 
\qquad c(u)\in \mathbb{F}_{2}^m.
\end{displaymath}
This representation is called the \emph{algebraic normal form} of $F$ and its 
degree $d^{\circ}(F)$ the \emph{algebraic degree} of the function $F$.\\
Besides, viewed as a function from the field $\mathbb{F}_{2^m}$ to itself, $F$ 
has a unique representation as a univariate polynomial over $\mathbb{F}_{2^m}$ 
of degree smaller than $2^m$:
\begin{displaymath}
F(x)=\sum_{i=0}^{2^m-1}c_ix^i,\quad c_i\in \mathbb{F}_{2^m}. 
\end{displaymath}
For any $k$, $0\le k\le 2^m-1$, the number $w_2(k)$ of the nonzero coefficients 
$k_s\in \{0,1\}$ in the binary expansion $\sum_{s=0}^{m-1}2^sk_s$ of $k$ is 
called the $2$-weight of $k$. 
The algebraic degree of $F$  is equal to the maximum 2-weight of the exponents 
$i$ of the polynomial $F(x)$ such that $c_i\neq 0$, that is 
$d^{\circ}(F)=\max_{0\le i \le n-1,c_i\ne 0}w_2(i)$ (see \cite{CC-carchzin}).

A function $F:\mathbb{F}_{2}^m\to \mathbb{F}_{2}^m$ is \emph{linear} if and only 
if $F(x)$ is a linearized polynomial over $\mathbb{F}_{2^m}$, that is,
\begin{displaymath}
\sum_{i=0}^{m-1}c_ix^{2^i},\quad c_i\in \mathbb{F}_{2^m}.
\end{displaymath}
The sum of a linear function and a constant is called an \emph{affine function}.

Let $F$ be a function from $\mathbb{F}_{2}^m$ to itself and $A_1$,
$A_2:\mathbb{F}_{2}^m\to \mathbb{F}_{2}^m$ be affine permutations. Then the 
functions $F$ and $A_1\circ F\circ A_2$ are called
\emph{affine equivalent}. Affine equivalent functions have the same algebraic 
degree (i.e. the algebraic degree is \emph{affine invariant}).

We shall say that the functions $F$ and $F^\prime$ are \emph{extended affine 
equivalent} (EA-equivalent) if $F^\prime=A_1\circ F\circ A_2+A$ for some affine 
permutations   $A_1$,
$A_2$ and an affine function $A$.
If $F$ is not affine, then  $F$ and $F^\prime$ have again the same algebraic 
degree.

For a function $F:\mathbb{F}_{2}^m\to \mathbb{F}_{2}^m$ and any elements $a,b\in 
\mathbb{F}_{2}^m$ we denote
$$
\delta_F(a,b)  =  |\{x\in \mathbb{F}_{2}^m: F(x+a)+F(x)=b\}|
$$
and
$$
\Delta_F  =  \{ \delta_F(a,b): a,b\in \mathbb{F}_{2}^m, a\ne 0\}.
$$
$F$ is called a \emph{differentially $\delta$-uniform} function if 
$\max_{a\in \mathbb{F}_{2}^{m*}, b\in \mathbb{F}_{2}^{m}}  \delta_F(a,b) \le 
\delta,$ where $\mathbb{F}_{2}^{m*}=\mathbb{F}_{2}^{m}\setminus \{0\}$.
For any $a,b\in \mathbb{F}_{2}^m$, the number $\delta_F(a,b)$ is even since if 
$x_0$ is a solution of the equation $F(x+a)+F(x)=b$ then $x_0+a$ is a solution 
too. Hence, $\delta\ge 2$. Differentially 2-uniform mappings are called 
\emph{almost perfect nonlinear}.

For any function $F:\mathbb{F}_{2}^m\to \mathbb{F}_{2}^m$, the values 
$$
\lambda_F(a,b) = \sum_{x\in \mathbb{F}_{2}^m}(-1)^{b\cdot F(x) + a\cdot 
x},\qquad a,b\in \mathbb{F}_{2}^m,
$$
do not depend on a particular choice of the inner product $"\cdot"$ in 
$\mathbb{F}_{2}^m$. If we identify $\mathbb{F}_{2}^m$ with $\mathbb{F}_{2^m}$ 
then we can take $x\cdot y=tr(xy)$, where $tr(x)=x+x^2+x^4+...+x^{2^{m-1}}$ is 
the trace function from $\mathbb{F}_{2^m}$ into $\mathbb{F}_{2}$.  Later we 
shall also need for a divisor $n$ of $m$ the relative trace 
$tr_{m/n}(x)=x+x^{2^n}+x^{2^{2n}}+...+x^{2^{(m/n-1)n}}$ and we denote 
$tr_{n}(x)=x+x^2+x^{2^2}+...+x^{2^{n-1}}$.
The set $\Lambda_F  =  \{ \lambda_F(a,b): a,b\in \mathbb{F}_{2}^m, b\ne 0\}$ is 
called the \emph{Walsh spectrum} of $F$ and the value 
$${\cal NL}(F)=2^{m-1}-\frac{1}{2} \max_{a\in \mathbb{F}_{2}^{m}, b\in 
\mathbb{F}_{2}^{m*}}  |\lambda_F(a,b) |$$  the \emph{nonlinearity} of the 
function $F$. 
The nonlinearity of any function $F$ satisfies the inequality 
$${\cal NL}(F)\le 2^{m-1}-2^\frac{m-1}{2}$$
(\cite{CV,sid}) and in case of equality $F$ is called \emph{almost bent}  or 
\emph{maximum nonlinear}.   For any AB function $F$, the Walsh spectrum 
$\Lambda_F$ equals $\{0,\pm 2^\frac{m+1}{2}\}$ as it is proven in \cite{CV}.

For EA-equivalent functions $F$ and $F^\prime$, we have 
$\Delta_F=\Delta_{F^\prime}$,  $\Lambda_F=\Lambda_{F^\prime}$ and if $F$ is a 
permutation then $\Delta_F=\Delta_{F^{-1}}$, $\Lambda_F=\Lambda_{F^{-1}}$ (see 
\cite{CC-carchzin}). Therefore, if $F$ is APN (resp. AB) and $F^\prime$ is 
EA-equivalent to either $F$ or $F^{-1}$ (if $F$ is a permutation), then 
$F^\prime$ is also APN (resp. AB).

Table 1 (resp. Table 2) gives all known values of exponents $d$ (up to 
EA-equivalence and up to taking the inverse when a function is a permutation) 
such that the power function $x^d$ is APN (resp. AB). 
\begin{center}
\small{Table 1

Known APN power functions  $x^d$ on $\mathbb{F}_{2^m}$.}

\begin{tabular}{|c|c|c|c|c|}
\hline
  & \footnotesize{Exponents $d$} & \footnotesize{Conditions} 
&\footnotesize{$w_2(d)$} & \footnotesize{Proven in}\\
\hline
\hline
\footnotesize{Gold functions} & \footnotesize{$2^i+1$} & 
\footnotesize{$\gcd(i,m)=1$} & \footnotesize{2}& 
\footnotesize{\cite{Gold,nyb}}\\

\hline
\footnotesize{Kasami functions} & \footnotesize{$2^{2i}-2^i+1$} & 
\footnotesize{$\gcd(i,m)=1$}  & \footnotesize{$i+1$}& 
\footnotesize{\cite{J-W,Kasami}}\\

\hline
\footnotesize{Welch function}  & \footnotesize{$2^t +3$ }&  
\footnotesize{$m=2t+1$}  &\footnotesize{3} &\footnotesize{\cite{D2}}\\
\hline
\footnotesize{Niho function}  &\footnotesize{$2^t+2^\frac{t}{2}-1$, $t$ even} & 
\footnotesize{$m=2t+1$}&\footnotesize{$(t+2)/2$} & \footnotesize{\cite{D1}}\\
 & \footnotesize{$2^t+2^\frac{3t+1}{2}-1$, $t$ odd} & & \footnotesize{$t+1$} & 
\\
\hline
\footnotesize{Inverse function } &\footnotesize{$2^{2t}-1 $}& 
\footnotesize{$m=2t+1$}&\footnotesize{$m-1$}& \footnotesize{\cite{B-D,nyb}}\\
\hline
\footnotesize{Dobbertin function}  & 
\footnotesize{$2^{4i}+2^{3i}+2^{2i}+2^{i}-1$} & 
\footnotesize{$m=5i$}&\footnotesize{$i+3$} & \footnotesize{\cite{D3}}\\
\hline
\end{tabular}
\end{center}

\begin{center}
\small{Table 2

Known AB power functions  $x^d$ on $\mathbb{F}_{2^m}$, $m$ odd.}

\begin{tabular}{|c|c|c|c|}
\hline
  & \footnotesize{Exponents $d$} & \footnotesize{Conditions} & 
\footnotesize{Proven in}\\
\hline
\hline
\footnotesize{Gold functions} & \footnotesize{$2^i+1$} & 
\footnotesize{$\gcd(i,m)=1$} & \footnotesize{\cite{Gold,nyb}}\\

\hline
\footnotesize{Kasami functions} & \footnotesize{$2^{2i}-2^i+1$} & 
\footnotesize{$\gcd(i,m)=1$}  & \footnotesize{\cite{Kasami}}\\

\hline
\footnotesize{Welch function}  & \footnotesize{$2^t +3$ }&  
\footnotesize{$m=2t+1$}  &\footnotesize{\cite{CCD2,CCD}}\\
\hline
\footnotesize{Niho function}  &\small{$2^t+2^\frac{t}{2}-1$, $t$ even} & 
\footnotesize{$m=2t+1$} & \footnotesize{\cite{H-X}}\\
 & \footnotesize{$2^t+2^\frac{3t+1}{2}-1$, $t$ odd} &  & \\
\hline
\end{tabular}\\
\end{center}
Every AB function is APN \cite{CV}. The converse is not true since AB functions 
exist only when  $m$ is odd while APN functions exist for $m$ even too. Besides, 
in the $m$ odd case, the Dobbertin APN function is not AB as proven in 
\cite{CCD}. Also the inverse APN function is not AB since it has the algebraic 
degree $m-1$ while the algebraic degree of any AB function is not greater than 
$(m+1)/2$ (see \cite{CC-carchzin}). When $m$ is even,  the inverse function 
$x^{2^m-2}$ is a differentially 4-uniform permutation  \cite{nyb} and has  the 
best known nonlinearity \cite{LW}, that is $2^{m-1}-2^\frac{m}{2}$ (see 
\cite{CCD,D4}). When $m\equiv 2$ [mod 4] and $\gcd (i,m)=2$, the functions 
$x^{2^i+1}$ and $x^{2^{2i}-2^i+1}$ also have the best known nonlinearity 
\cite{Gold,Kasami} and when  $\gcd (i,m)=s$ and $m/s$ is odd, these functions 
have three valued Walsh spectrum $\{0,\pm 2^\frac{m+s}{2}\}$ (\cite{nyb,T}). 

\section{Carlet-Charpin-Zinoviev equivalence of functions}
For a function $F$ from $\mathbb{F}_{2}^m$ to itself, we denote by $G_F$  the 
\emph{ graph of the function} $F$:
\begin{displaymath}
G_F=\{(x,F(x)):x\in\mathbb{F}_{2}^m\}\subset \mathbb{F}_{2}^{2m}.
\end{displaymath}
The property of stability of APN and AB functions given in Proposition~3 of 
\cite{CC-carchzin}  leads to the definition of the following equivalence 
relation of functions.
\newtheorem{definition}{Definition}
\begin{definition} \emph{(\textbf{CCZ-equivalence})} We say that  functions 
$F,F^\prime:\mathbb{F}_{2}^m\to \mathbb{F}_{2}^m$ are 
\emph{Carlet-Charpin-Zinoviev  equivalent} if there exists a linear (affine) 
permutation ${\cal L}:\mathbb{F}_{2}^{2m}\to \mathbb{F}_{2}^{2m}$ such that 
${\cal L}(G_F)=G_{F^\prime}$.
\end{definition}
We shall consider only the case of linear functions, but all proofs for the 
statements related to the CCZ-equivalence can be easily extended for the case of 
affine functions too.

A linear function ${\cal L}:\mathbb{F}_{2}^{2m}\to \mathbb{F}_{2}^{2m}$  can be 
considered as a pair of linear functions $L_1,L_2:\mathbb{F}_{2}^{2m}\to 
\mathbb{F}_{2}^m$ such that ${\cal L}(x,y)=(L_1(x,y),L_2(x,y))$ for all 
$x,y\in\mathbb{F}_{2}^m$. Then for a function $F:\mathbb{F}_{2}^m\to 
\mathbb{F}_{2}^m$ we have ${\cal L}(x,F(x))=(F_1(x),F_2(x))$, where
\begin{equation}\label{F1}
F_1(x)=L_1(x,F(x)), 
\end{equation}
\begin{equation}\label{F2}
F_2(x)=L_2(x,F(x)). 
\end{equation}
Hence, the set ${\cal L}(G_F)=\{(F_1(x),F_2(x)):x\in\mathbb{F}_{2}^m\}$ is the
graph of a function $F^\prime$ if and only if the function $F_1$ is a 
permutation. If ${\cal L}$ and $F_1$ are permutations then ${\cal 
L}(G_F)=G_{F^\prime}$, where $F^\prime=F_2\circ F_1^{-1}$, and the functions $F$ 
and $F^\prime$ are CCZ-equivalent. \\

In the proposition below, we give a slightly different approach to the 
CCZ-equivalence, that we shall use in the constructions of APN polynomials which 
will be EA-inequivalent to power functions. Recall that a set 
$G\subset\mathbb{F}_{2}^{2m}$ is {\em transversal} to a subgroup $V$ of 
$\mathbb{F}_{2}^{2m}$ if $|G\cap (u+~V)|=~1$ for any $u\in \mathbb{F}_{2}^{2m}$.
\newtheorem{Proposition}{Proposition}
\begin{Proposition} \label{prop2}
Let $F:\mathbb{F}_{2}^m\to\mathbb{F}_{2}^m$ and  $G\subset\mathbb{F}_{2}^{2m}$. 
Then the set $G$ is the graph of a function CCZ-equivalent to $F$ if and only if 
there exists a linear permutation  ${\cal L} :\mathbb{F}_{2}^{2m}\to 
\mathbb{F}_{2}^{2m}$ such that $G={\cal L} (G_F)$ and  $G_{F}$ is transversal to 
$V^\prime={\cal L}^{-1} (V)$, where $V=\{(0,x):x\in\mathbb{F}_{2}^m\}.$
\end{Proposition}
\emph{Proof.} The condition that there exists a linear permutation ${\cal L} 
:\mathbb{F}_{2}^{2m}\to \mathbb{F}_{2}^{2m}$ such that $G={\cal L} (G_F)$ is 
clearly necessary.
Let us denote $U=\{(x,0):x\in\mathbb{F}_{2}^m\}$, 
$V=\{(0,x):x\in\mathbb{F}_{2}^m\}$, ${\cal L}^{-1}(U)=U^\prime$  and ${\cal 
L}^{-1}(V)=V^\prime$. Then $G$ is the graph of a function if and only if $|G\cap 
(u+V)|=1$ for any $u\in U$; that is, if and only if  $|{\cal L}^{-1}(G)\cap 
(u^\prime+V^\prime)|=1$ for any $u^\prime\in U^\prime$. Hence, $G$ is the graph 
of a function CCZ-equivalent to $F$ if and only if $G_{F}$ is transversal to 
$V^\prime$.  \hfill $\Box$\\

Let  $F$ be an arbitrary function on $\mathbb{F}_{2}^m$ and $V^\prime$ be an 
arbitrary subgroup of $\mathbb{F}_{2}^{2m}$. If we denote
$$H_a=\{F(x+a)+F(x): x\in  \mathbb{F}_2^m\},\qquad a\in \mathbb{F}_2^{m*},$$
$$A_F=\{(a,F(x+a)+F(x)): a\in \mathbb{F}_{2^m}^{*}, x\in \mathbb{F}_{2^m}\},$$
then $A_F=\bigcup_{a\in \mathbb{F}_{2^m}^{*}} (a,H_a).$ It is  easy to note that 
$A_F\cap V^\prime=\emptyset$  if and only if $G_F$ is transversal to $V^\prime$.
We shall see that, to construct a function CCZ-equivalent but EA-inequivalent to 
the function $F$, it is necessary that the subgroup $V'$ in Proposition 
\ref{prop2} is different from $V$ and as soon as we have a subroup $V'$ such 
that $G_F$ is transversal to $V'$ we can find a linear permutation ${\cal L}$ 
that $V^\prime={\cal L}^{-1} (V)$. However, even if we have such a subgroup it 
is difficult to determine whether the resulting function is EA-inequivalent to 
$F$.\\

The property of stability of APN and AB mappings (see \cite{CC-carchzin}) can be 
easily generalized to all functions (not necessarily APN or AB) as follows:
\begin{Proposition}\label{prop1}
Let $F,F^\prime:\mathbb{F}_{2}^m\to \mathbb{F}_{2}^m$ be CCZ-equivalent 
functions. Then  $\Delta_F=\Delta_{F^\prime}$ and 
$\Lambda_F=\Lambda_{F^\prime}$. In particular,  $F$ is an APN (resp. AB) 
function if and only if $F^\prime$ is APN (resp. AB).
\end{Proposition}
\emph{Proof.} If $F,F^\prime:\mathbb{F}_{2}^m\to \mathbb{F}_{2}^m$ are 
CCZ-equivalent, then $F^\prime=F_2\circ F_1^{-1}$ for a certain linear 
permutation ${\cal L}=(L_1,L_2)$, where $F_1,F_2$ are defined by (\ref{F1}) and 
(\ref{F2}). For any $(a,b)\in \mathbb{F}_{2}^{2m}$ we have
{\setlength\arraycolsep{2pt}
\begin{eqnarray*}
\lambda_{F}(a,b) & = & \sum_{x\in \mathbb{F}_{2}^m}(-1)^{b\cdot F(x) + a\cdot 
x}=\sum_{x\in \mathbb{F}_{2}^m}(-1)^{(a,b)\cdot (x,F(x))} \\
 & = &  \sum_{x\in \mathbb{F}_{2}^m}(-1)^{(a,b)\cdot {\cal 
L}^{-1}(F_1(x),F_2(x))}
 =\sum_{x\in \mathbb{F}_{2}^m}(-1)^{{\cal L}^{-1*}(a,b)\cdot (x,F_2\circ 
F_1^{-1}(x))}\\
& = & \lambda_{F^\prime}({\cal L}^{-1*}(a,b)),
\end{eqnarray*}}
where ${\cal L}^{-1*}$ is the adjoint operator of ${\cal L}^{-1}$ (i.e. $x\cdot 
{\cal L}^{-1}(y)={\cal L}^{-1*}(x)\cdot y$, for any $(x,y)\in 
\mathbb{F}_{2}^{2m}$~; if ``$\cdot$" is the usual inner product, then ${\cal 
L}^{-1*}$ is the linear permutation whose matrix is transposed of that of ${\cal 
L}^{-1}$).
Hence, $\Lambda_F=\Lambda_{F^\prime}$.\\
The proof that $\Delta_F=\Delta_{F^\prime}$ for arbitrary functions $F,F^\prime$ 
is completely the same as in the case when one of the functions $F,F^\prime$ is 
APN (see \cite{CC-carchzin}).  \hfill $\Box$\\

\noindent \textbf{Remark 1} Obviously, CCZ-equivalence can be defined for 
functions between any two groups $H_1$ and $H_2$. For a function $F:H_1\to H_2$ 
we can consider  the set of the values $\delta_F(a,b)=|\{x\in 
H_1:F(x+a)-F(x)=b\}|$, $a\in H_1\backslash\{0\}$, $b\in H_2$, and  if  the 
groups $H_1$ and $H_2$ are Abelian, then the discrete Fourier transform of $F$ 
can also be defined. In this more general case CCZ-equivalent functions again 
have the same differential and linear properties. One can find results related 
to nonlinear functions on Abelian groups in \cite{CD,P}.
\hfill $\Diamond$\\

Since CCZ-equivalent functions have the same differential uniformity and the 
same nonlinearity, the resistance of a function to linear and differential 
attacks is CCZ-invariant. CCZ-equivalent functions have also the same 
weakness/strength with respect to algebraic attacks. Indeed, let functions 
$F,F^\prime:\mathbb{F}_{2}^m\to \mathbb{F}_{2}^m$ be CCZ-equivalent. Then 
$F^\prime=F_2\circ F_1^{-1}$, where $F_1,F_2$ are defined by (\ref{F1}) and 
(\ref{F2}) for a certain linear permutation ${\cal L}=(L_1,L_2)$.  If there 
exists a nonzero function $\psi:\mathbb{F}_{2}^{2m}\to \mathbb{F}_{2}$ of low 
degree such that 
\begin{equation}\label{ajout1}
\psi(x,F(x))=0, \quad\forall x\in\mathbb{F}_{2}^{m},
\end{equation}
then $\psi$ could be used in an algebraic attack \cite{CourtPier}.  The function 
$\psi\circ {\cal L}^{-1}$ has the same degree as $\psi$ and Relation 
(\ref{ajout1}) is equivalent to
\begin{displaymath}
\psi \circ {\cal L}^{-1}(F_1(x),F_2(x))=0,\quad\forall x\in\mathbb{F}_{2}^{m}, 
\end{displaymath}
which implies
\begin{displaymath}
\psi \circ {\cal L}^{-1}(x,F^\prime(x))=0,\quad\forall x\in\mathbb{F}_{2}^{m}
\end{displaymath}
by replacing $x$ by $F_{1}^{-1}(x)$.
Hence, $\psi \circ {\cal L}^{-1}$ could be used in an algebraic attack on 
$F^\prime$, and {\em vice versa}. Therefore, the resistance of a function to 
algebraic attacks is also CCZ-invariant.\\
Since any permutation is CCZ-equivalent to its inverse then obviously the 
algebraic degree of a function is not CCZ-invariant (while it is EA-invariant as 
noted above). For example, if $F:\mathbb{F}_{2}^m\to \mathbb{F}_{2}^m$ is a Gold 
AB function  then 
$d^{\circ}(F)=2$ and $d^{\circ}(F^{-1})=\frac{m+1}{2}$ as proven in 
\cite{nyb}.\\
 
EA-equivalent functions are CCZ-equivalent and if a function $F$ is a 
permutation then $F$ is CCZ-equivalent to $F^{-1}$ \cite{CC-carchzin}. More 
precisely:
\begin{Proposition}\label{prop0} Let 
$F,F^\prime:\mathbb{F}_2^m\to\mathbb{F}_2^m$. The  function $F^\prime$ is 
EA-equivalent to the function $F$ or to the inverse of $F$ (if it exists) if and 
only if there exists a linear permutation ${\cal L}=(L_1,L_2)$ on 
$\mathbb{F}_2^{2m}$ such that ${\cal L}(G_F)=G_{F^\prime}$ and the function 
$L_1$ depends only on one variable, i.e. $L_1(x,y)=L(x)$ or $L_1(x,y)=L(y)$.
\end{Proposition}
\emph{Proof.}  Let ${\cal L}(G_F)=G_{F^\prime}$ for some linear permutation 
${\cal L}=(L_1,L_2):\mathbb{F}_2^{2m}\to\mathbb{F}_2^{2m}$ and $L_1(x,y)=L(y)$, 
$L_2(x,y)=R_1(x)+R_2(y)$, where $L$, $R_1$, $R_2:\mathbb{F}_2^m\to 
\mathbb{F}_2^m$ are linear. Then 
$$F_1(x)=L_1(x,F(x))=L\circ F (x), $$ 
$$F_2(x)=L_2(x,F(x))=R_1(x)+R_2\circ F(x).$$
The function $F_1$ is a permutation, since ${\cal L}(G_F)$ is the graph of a 
function. Therefore, $L$ and $F$ have to be permutations.  
On the other hand, the system
\begin{displaymath}
\left\{ \begin{array}{lll}
L(y)&=&0\\
R_1(x)+R_2(y)&=&0
\end{array}\right.
\end{displaymath}
has only $(0,0)$ solution, since ${\cal L}=(L_1,L_2)$ is a permutation. But $L$ 
is also a permutation. Therefore, the linear function $R_1$ has to be a 
permutation too.\\
We have
{\setlength\arraycolsep{2pt}
\begin{eqnarray*}
F^\prime (x)&=&F_2\circ F_1^{-1}(x)=R_1\circ F^{-1}\circ L^{-1}(x)+R_2\circ 
F\circ F^{-1}\circ L^{-1}(x)\\
&=&R_1\circ F^{-1}\circ L^{-1}(x)+R_2\circ L^{-1}(x).
\end{eqnarray*}}
Thus, $F^\prime$ is EA-equivalent to $F^{-1}$.\\
The proof that $F^\prime$ is EA-equivalent to $F$, when $L_1(x,y)=L(x)$, is 
similar.

Conversely, let $F^\prime=R_1\circ F\circ R_2+R$ or $F^\prime=R_1\circ 
F^{-1}\circ R_2+R$ for some linear permutations $R_1,R_2$ and for some linear 
function $R$. Then take ${\cal L}(x,y)=(R_2^{-1}(x),R_1(y)+R\circ R_2^{-1}(x))$ 
in the first case and in the second case take ${\cal L}(x,y)=(R_2^{-1}(y),R\circ 
R_2^{-1}(y)+R_1(x))$. Obviously, all conditions are satisfied.
\hfill $\Box$\\

\noindent \textbf{Remark 2} Proposition \ref{prop0} shows that, for a function 
$F:\mathbb{F}_{2}^m\to \mathbb{F}_{2}^m$, if ${\cal L}=(L_1,L_2)$ and ${\cal 
L}^\prime=(L_1,L_2^\prime)$ are linear permutations on $\mathbb{F}_2^{2m}$ such 
that the function $L_1(x,F(x))$ is a permutation, then the functions defined by 
the graphs ${\cal L}(G_F)$ and ${\cal L}^\prime(G_F)$ are EA-equivalent.
Indeed, we have
$({\cal L}^\prime\circ {\cal L}^{-1})\circ {\cal L}(G_F)={\cal L}^\prime(G_F) $ 
and, denoting  $({\cal L}^\prime\circ {\cal L}^{-1})=(S_1,S_2)$, we have 
$S_{1}(x,y)=x$. This last equality can be easily checked by considering the 
linear functions ${\cal L}$, ${\cal L}^\prime$  and ${\cal L}^{-1}$ as 
$(2m)\times (2m)$ matrices
\begin{displaymath}
{\cal L}=
\left ( \begin{array}{lll}
R_1 & R_2\\
T_1 & T_2
\end{array}\right),\quad
{\cal L}^\prime=
\left ( \begin{array}{lll}
R_1 & R_2\\
T_1^\prime & T_2^\prime
\end{array}\right),\quad
{\cal L}^{-1}=
\left ( \begin{array}{lll}
A_1 & A_2\\
A_3 & A_4
\end{array}\right).
\end{displaymath}
We have
\begin{displaymath}
\left ( \begin{array}{lll}
R_1 & R_2\\
T_1 & T_2
\end{array}\right)\times
\left ( \begin{array}{lll}
A_1 & A_2\\
A_3 & A_4
\end{array}\right)=
\left ( \begin{array}{lll}
R_1A_1 + R_2A_3 & R_1A_2 + R_2A_4\\
T_1A_1 + T_2A_3 & T_1A_2 + T_2A_4
\end{array}\right)=
\left ( \begin{array}{lll}
I & 0\\
0 & I
\end{array}\right),
\end{displaymath}where $I$ is the identity matrix and $0$ is the 0-matrix of 
order $m$, and this implies
\begin{displaymath}
\left ( \begin{array}{lll}
R_1 & R_2\\
T_1^\prime & T_2^\prime
\end{array}\right)\times
\left ( \begin{array}{lll}
A_1 & A_2\\
A_3 & A_4
\end{array}\right)=
\left ( \begin{array}{ccc}
I & 0\\
T_1^\prime A_1 + T_2^\prime A_3 & T_1^\prime A_2 + T_2^\prime A_4
\end{array}\right).
\end{displaymath} 
\hfill $\Diamond$\\

Proposition \ref{prop0} shows that if we want to construct  functions $F^\prime$ 
which are CCZ-equivalent to a function $F$ and EA-inequivalent to  both $F$ and 
$F^{-1}$ (if $F^{-1}$ exists), then we have to use a linear function $L_1(x,y)$ 
depending on both variables. However, this condition is  not sufficient as the 
following example shows.\\
 
\noindent \textbf{Example 1 }  Let $m=2n+1$ and $s\equiv n$ [mod $2]$. Then the 
linear function 
$${\cal L}(x,y)=(x+tr(x)+\sum_{j=0}^{n-s}y^{2^{2j+s}},y+tr(x))$$ 
is a permutation on $\mathbb{F}_{2^m}^2$ since the kernel of ${\cal L}$ is 
$\{(0,0)\}$ (this can be checked by considering the cases $tr(x)=0$ and 
$tr(x)=1$). For the Gold AB function $x^3$ the function 
$$F_1(x)=x+tr(x)+\sum_{j=0}^{n-s}(x^3)^{2^{2j+s}}$$
is a permutation on $\mathbb{F}_{2^m}$. Indeed, denoting 
$L(y)=\sum_{j=0}^{n-s}y^{2^{2j+s}}$ we have $L(y+y^2)=y+tr(y)$ and 
$L((y+1)^3)=L(y^3)+y+tr(y)+1$ since $L(1)=1$. Thus $F_1(x)=L((x+1)^3)+1$ and 
$F_1$ is a permutation if $L$ is bijective. The equation $z=L(y)$ implies 
$z+z^2=y+tr(y)$ and $tr(z)=tr(y)$. Therefore, $L$ is a permutation and 
$L^{-1}(x)=x+x^2+tr(x)$, $F_1^{-1}(x)=[L^{-1}(x+1)]^\frac{1}{3}+1$. Finally, we 
get
{\setlength\arraycolsep{2pt}
\begin{eqnarray*}
F^\prime(x)&=&F_2\circ 
F_1^{-1}(x)=([L^{-1}(x+1)]^\frac{1}{3}+1)^3+tr([L^{-1}(x+1)]^\frac{1}{3}+1)\\
&=&L^{-1}(x+1)+[L^{-1}(x+1)]^\frac{2}{3}+[L^{-1}(x+1)]^\frac{1}{3}+tr([L^{-1}(x+1)]^\frac{1}{3})\\
&=&L^{-1}(x+1)+L^{-1}([L^{-1}(x+1)]^\frac{1}{3}).
\end{eqnarray*}}
Thus, both functions $L_1$ and $L_2$ depend on two variables, but the function 
$F^\prime$ is  EA-equivalent to the inverse of $x^3$.\\
This example can be generalized for any Gold AB function by replacing 
$L^{-1}(x)=x+x^{2^i}+tr(x)$ and $x^3$ by $x^{2^i+1}$.
\hfill $\Diamond$\\

If we want to classify all functions  CCZ-equivalent to a given one
$F$, then we should characterize all permutations of the form $L\circ
F+L^\prime$, where $L,L^\prime$ are linear. Indeed, we need to know
which linear mapping $L_{1}~:\mathbb{F}_{2}^{2m}\to \mathbb{F}_{2}^m$
is such that the function $F_1(x)=L_1(x,F(x))$ is a permutation. Clearly, 
$L_{1}$ can be written uniquely in the form $L_{1}(x,y)=L(y)+L^\prime(x)$.
If $F_1$ is a permutation then there exists a linear function $L_2(x,y)$ such 
that the linear function $(L_1,L_2)(x,y)$ is a permutation too. Indeed, 
$L_{1}(x,F(x))$ being a permutation, $L_1$ is onto and the kernel of $L_1$ has 
then dimension $m$. We can take for $L_2$ any linear permutation between $Ker 
(L_1)$ and $\mathbb{F}_{2}^m$, extended to $\mathbb{F}_{2}^{2m}$ by 
$L_2(x+y)=L_2(x)$ for all $x\in Ker (L_1)$, $y\in E$, where $E$ is an 
$m$-dimensional subspace of $F_2^{2m}$ such that $E\oplus Ker 
(L_1)=\mathbb{F}_{2}^{2m}$. Conversely, any linear function $L_2$  such that 
$(L_1,L_2)$ is a permutation has this form. Indeed, $L_{2}$ being onto, it has 
also an $m$-dimensional kernel, and $(L_{1},L_{2})$ being one to one, the 
kernels of $L_{1}$ and $L_{2}$ have trivial intersection. This proves that $Ker 
(L_1)\oplus Ker (L_2)=\mathbb{F}_{2}^{2m}$ and that we can take $E=Ker (L_2)$.\\

The following proposition gives a sufficient condition for a function to be 
EA-inequivalent to power functions.
\begin{Proposition}\label{prop5}
Let $F$ be a function from $\mathbb{F}_{2^m}$ to itself. If there exists an 
element $c\in\mathbb{F}_{2^m}^*$ such that $d^{\circ}(tr(cF))\notin 
\{0,1,d^{\circ}(F)\}$, then
 $F$ is EA-inequivalent to power functions. 
\end{Proposition}
\emph{Proof.} It is proven in \cite{Cthesis} that for any power function $x^d$ 
and for any $c\in \mathbb{F}_{2^m}$, either the function $tr(cx^d)$ completely 
vanishes or it has exactly the algebraic degree $w_2(d)$. 
Thus, for any function $F$ which is affine equivalent to a power 
function, we have $d^{\circ}(tr(cF))\in \{0,d^{\circ}(F)\}$, $c\in 
\mathbb{F}_{2^m}^*$. Therefore, if $F$ is EA-equivalent to a power function then 
$d^{\circ}(tr(cF))\in \{0,1,d^{\circ}(F)\}$, for every   $c\in 
\mathbb{F}_{2^m}^*$. \hfill $\Box$

\section{CCZ-equivalence and the Gold functions}
In the propositions below we give a characterization of permutations  $L\circ
F+L^\prime$ when $F$ is a Gold function. This characterization is not complete 
but it is useful for constructions of new APN and AB polynomials.
\begin{Proposition}\label{prop3}
Let $F:\mathbb{F}_{2^m}\to \mathbb{F}_{2^m}$, $F(x)=L(x^{2^i+1})+L^\prime(x)$, 
where  $L,L^\prime$ are linear and $\gcd(m,i)=1$. Then $F$ is a permutation if 
and only if, for every $u\neq 0$ in $\mathbb{F}_{2^m}$ and every $v$ such that 
$tr(v)=tr(1)$, the condition $L(u^{2^i+1}v)\neq L'(u)$ holds. 
\end{Proposition}
\emph{Proof.}
For any $u\in\mathbb{F}_{2^m}^*$ we have
{\setlength\arraycolsep{2pt}
\begin{eqnarray*}
&&F(x)+F(x+u)=L(x^{2^i+1})+L^\prime(x)+L((x+u)^{2^i+1})+L^\prime(x+u)\\
&&=L\left(ux^{2^i}+u^{2^i}x+ 
u^{2^i+1}\right)+L'(u)=L\left(u^{2^i+1}\left((x/u)^{2^i}+x/u+1\right)\right)+L'(u).
\end{eqnarray*}}
When $x$ ranges over $\mathbb{F}_{2^m}$ then $(x/u)^{2^i}+x/u+1$ ranges over  
$\{v\in \mathbb{F}_{2^m}: tr(v)=tr(1)\}$, since $\gcd(m,i)=1$.  Hence, $F$ is 
injective (i.e. is a permutation) if and only if $L(u^{2^i+1}v)\neq L'(u)$ for 
every  $u\neq 0$ and every $v$ such that $tr(v)=tr(1)$.\hfill $\Box$\\
 
\noindent \textbf{Remark 3 } If $m$ is even, then, up to affine equivalence and 
without loss of generality, we can consider only functions of the type 
$L(x^{2^i+1})+x$. Indeed, if  $F(x)=L(x^{2^i+1})+L^\prime(x)$ is a permutation 
on $\mathbb{F}_{2^m}$ and $m$ is even, then it follows from Proposition 
\ref{prop3} that $L'$ must be a permutation (take $v=0$). Then the function 
$F^\prime(x)=L^{\prime -1}\circ F(x)=L^{\prime -1}\circ L(x^{2^i+1})+x$ is a 
permutation if and only if $F$ is a permutation. Moreover if a function 
$L_{2}~:\mathbb{F}_{2^m}^{2}\to \mathbb{F}_{2^m}$ is such that the function 
$(L(y)+L^\prime(x),L_2(x,y))$ is a permutation on $\mathbb{F}_{2^m}^{2}$, then 
obviously $(L^{\prime -1}\circ L(y)+x,L_2(x,y))$ is a permutation too; note also 
that $(L(y)+x,y)$ is  a permutation for any linear function $L$. Thus the 
function $x^{2^i+1}$ is CCZ-equivalent to the functions $F_2\circ F^{\prime -1}$ 
and  $F_2\circ F^{-1}$, where $F_2(x)=L_2(x,x^{2^i+1})$.  We have $F_2\circ 
F^{\prime -1}(x)=F_2\circ F^{-1}\circ L^\prime(x).$ Therefore, $F_2\circ 
F^{\prime -1}$ and  $F_2\circ F^{-1}$ are affine equivalent. \hfill $\Diamond$
\begin{Proposition}\label{prop4}
Let $F:\mathbb{F}_{2^m}\to \mathbb{F}_{2^m}$, $F(x)=L(x^{2^i+1})+x$, where $L$ 
is linear, $m$ even and $\gcd(m,i)=1$. Let $L^*$ be the adjoint operator of $L$. 
Then $F$ is a permutation if and only if, for every $v\in \mathbb{F}_{2^m}$ such 
that $L^*(v)\neq 0$, there exists $u\in \mathbb{F}_{2^m}$ such that 
$L^*(v)=u^{2^{i}+1}$ and $tr_{m/2}(\frac{v}{u})\neq 0$, where $tr_{m/2}$ is the 
trace function from $\mathbb{F}_{2^m}$ to $\mathbb{F}_{2^{2}}$. 
\end{Proposition}
\emph{Proof.} 
The function $F$ is a  permutation if and only if, for every $v\neq 0$, the 
Boolean function  $tr(v(L(x^{2^i+1})+x))$ is balanced (see e.g. \cite{Cbook}). 
Let $L^*$ be the adjoint operator of $L$, then we have
$$tr(v(L(x^{2^i+1})+x))=tr(L^*(v)x^{2^i+1}+vx).$$
If $L^*(v)= 0$, then the function $tr(v(L(x^{2^i+1})+x))$ is balanced. If 
$L^*(v)\neq 0$, the quadratic function $tr(L^*(v)x^{2^i+1}+vx)$ has associated 
symplectic form~:
$$\varphi(x,y)=tr(L^*(v)x^{2^i}y+L^*(v)xy^{2^{i}})=tr((L^*(v)x^{2^i}+(L^*(v)x)^{2^{m-i}})y)~,$$
which has kernel~:
$${\cal E}=\{x\in \mathbb{F}_{2^m}: L^*(v)x^{2^i}+(L^*(v)x)^{2^{m-i}}=0\}=$$
$$=\{x\in \mathbb{F}_{2^m}: L^*(v)^{2^{i}}x^{2^{2i}}+L^*(v)x=0\}=\{0\}\cup 
\{x\in \mathbb{F}_{2^m}:L^*(v)^{2^i-1}x^{2^{2i}-1}=1\}.$$
A quadratic function is balanced if and only if its restriction to the kernel of 
its associated symplectic form is not constant \cite{Cthesis,Cbook1}. 
The restriction of $tr(L^*(v)x^{2^i+1})$ to ${\cal E}$ is null. Indeed, 
 $L^*(v)^{2^i-1}x^{2^{2i}-1}=1$ implies that the order of $L^*(v)x^{2^i+1}$ 
 divides $2^i-1$ and since $\gcd (i,m)=1$ then $L^*(v)x^{2^i+1}\in \mathbb{F}_2$ 
 and the trace function is null on $\mathbb{F}_{2}$, since $m$ is even.
Hence, the function $L(x^{2^i+1})+x$ is a  permutation if and only if every $v$ 
such that $L^*(v)\neq 0$ lies outside the dual of $\{0\}\cup \{x\in 
\mathbb{F}_{2^m}:L^*(v)^{2^i-1}x^{2^{2i}-1}=1\}$. Equivalently, the function 
$L(x^{2^i+1})+x$ is a  permutation if and only if, for every $v$ such that 
$L^*(v)\neq 0$ the following two conditions hold:\\
1) $L^*(v)^{2^i-1}$ belongs to $\{x^{2^{2i}-1}~: x\in \mathbb{F}_{2^m}\}$ 
(otherwise, ${\cal E}$ would be trivial), say $L^*(v)^{2^i-1}=u^{2^{2i}-1}$, 
i.e. $L^*(v)=u^{2^{i}+1}$ (since the mapping $y\rightarrow y^{2^{i}-1}$ is a 
permutation); in this case 
${\cal E}=\{0\}\cup\{x\in \mathbb{F}_{2^m}:(ux)^{2^{2i}-1}=1\}= 
\frac{1}{u}\{y\in \mathbb{F}_{2^m}:y^{2^{2i}}=y\}=\frac{1}{u}\mathbb{F}_{2^j}$,
where $j=\gcd(2i,m)=2$, hence ${\cal E}=\frac{1}{u}\mathbb{F}_{4}$; \\
2) $v$ lies outside the dual of ${\cal E}$, that is, $L(vx)\ne 0$ for some 
 $x\in{\cal E}$ .\\ 
 For every $z\in \mathbb{F}_{2^m}$ and every $y\in \mathbb{F}_{2^{2}}$ we have 
 $$tr(z\frac{1}{u}y)=tr_{2}(tr_{m/2}(\frac{z}{u}y))  =tr_{2}(y\, 
 tr_{m/2}(\frac{z}{u})).$$ 
 Hence, the function 
$L(x^{2^i+1})+x$ is a  permutation if and only if every $v$ such that 
$L^*(v)\neq 0$ satisfies $L^*(v)=u^{2^{i}+1}$ for some $u$ and 
$tr_{m/2}(\frac{v}{u})\neq 0$.  \hfill $\Box$\\

\section{New cases of AB and APN functions}
The next theorems show that  CCZ-equivalent functions are not necessarily 
EA-equivalent (including the equivalence to the inverse).  They lead to infinite 
classes of almost bent and almost perfect nonlinear polynomials, which are 
EA-inequivalent to power functions. 

For the function $F^\prime:\mathbb{F}_{2^m}\to \mathbb{F}_{2^m}$, 
$F(x)=x^{2^i+1}$, $\gcd(m,i)=1$, and for any $a\in \mathbb{F}_{2^m}^{*}$ the set
$$H_a=\{F(x+a)+F(x): x\in  \mathbb{F}_{2^m}\}=
\left\{ \begin{array}{ll}
\{x\in  \mathbb{F}_{2^m}: tr(a^{-(2^i+1)}x)=1\}& \quad\textrm{if } m \textrm{ is 
odd}\\
\{x\in  \mathbb{F}_{2^m}: tr(a^{-(2^i+1)}x)=0\}& \quad\textrm{if } m \textrm{ is 
even}
\end{array}\right.
$$
is a linear hyperplane when  $m$ is even and a complement of a linear hyperplane 
when  $m$ is  odd (see  \cite{BFF} or Proposition \ref{prop3}).\\
We can use Proposition \ref{prop2}. For any $a\in \mathbb{F}_{2^m}^{*}$ the set 
\begin{displaymath}
V^\prime=
\left\{ \begin{array}{ll}
(0,\mathbb{F}_{2^m}\backslash H_a)\cup (a,\mathbb{F}_{2^m}\backslash H_a)& 
\quad\textrm{if } m \textrm{ is odd}\\ 
(0,H_a)\cup (a,\mathbb{F}_{2^m}\backslash H_a)& \quad\textrm{if } m \textrm{ is 
even}
\end{array}\right.
 \end{displaymath}
is a  subgroup of $\mathbb{F}_{2^m}^2$ such that $G_F$ is transversal to 
$V^\prime$ since  $A_F\cap V^\prime=\emptyset$.
\newtheorem{Theorem}{Theorem}
\begin{Theorem}\label{thm1}
The function $F^\prime:\mathbb{F}_{2^m}\to \mathbb{F}_{2^m}$, $m>3$ odd,
\begin{displaymath}
F^\prime(x)=x^{2^i+1}+(x^{2^i}+x)tr(x^{2^i+1}+x), \textbf{  } \gcd(m,i)=1,
\end{displaymath}
is  an AB function which is EA-inequivalent to any power function.
\end{Theorem}
\emph{Proof.} It is easy (but lengthy) to prove that for 
$a\in\mathbb{F}_{2^m}^*$ the linear mapping
$${\cal L}(x,y)=(x+a\ tr(a^{-1}x)+a\ 
tr(a^{-(2^i+1)}y),y+a^{2^i+1}tr(a^{-(2^i+1)}y)+a^{2^i+1}tr(a^{-1}x))$$
is an involution.
We have
\begin{displaymath}
{\cal L}(0,y)=
\left\{ \begin{array}{ll}
(a,y+a^{2^i+1})& \textrm{if } \textrm{ } y\in H_a\\
(0,y)& \textrm{if } \textrm{ } y\in \mathbb{F}_{2^m}\backslash H_a
\end{array}\right..
\end{displaymath}
Since $tr(a^{-(2^i+1)}a^{2^i+1})=1$, then $a^{2^i+1}\in H_a$. The set $H_a$ is a 
complement of a linear hyperplane, hence the sum of any two elements from $H_a$ 
belongs to $\mathbb{F}_{2^m}\backslash H_a$. Therefore, ${\cal L}^{-1}(V)={\cal 
L}(V)=V^\prime$ and by Proposition \ref{prop2} the function $F$ is 
CCZ-equivalent to $F_2\circ F_1^{-1}$, where
{\setlength\arraycolsep{2pt}
\begin{eqnarray*}
F_1(x)&=&L_1(x,F(x))= x+a\ tr(x/a)+a\ tr((x/a)^{2^i+1}),\\
F_2(x)&=&L_2(x,F(x))= x^{2^i+1}+a^{2^i+1}tr((x/a)^{2^i+1})+a^{2^i+1}tr(x/a).
\end{eqnarray*}}
The function $F_1$ is an involution (this proves again that it is bijective):
{\setlength\arraycolsep{2pt}
\begin{eqnarray*}
F_1\circ F_1(x)&=&x+a\ tr(x/a)+a\ tr((x/a)^{2^i+1})+a\ tr(a^{-1}( x+a\ tr(x/a)\\
&+&a\ tr((x/a)^{2^i+1})))+a\ tr(a^{-(2^i+1)}( x+a\ tr(x/a)+a\ 
tr((x/a)^{2^i+1}))^{2^i+1})\\
&=&x+3a\ tr(x/a)+2a\ tr((x/a)^{2^i+1})+a\ 
tr(a^{-(2^i+1)}(x^{2^i+1}+(ax^{2^i}+a^{2^i}x\\
&+&a^{2^i+1})(tr(x/a)+tr((x/a)^{2^i+1}))))=x+a\ tr(x/a)+a\ tr((x/a)^{2^i+1})\\
&+&a\ tr((x/a)^{2^i}+(x/a)+1)(tr(x/a)+tr((x/a)^{2^i+1}))=x,
\end{eqnarray*}}
since
$$( x+a\ tr(x/a)+a\ tr((x/a)^{2^i+1}))^{2^i+1}=x^{2^i+1}+(ax^{2^i}+a^{2^i}x
+a^{2^i+1})(tr(x/a)+tr((x/a)^{2^i+1}))$$
$$+2a^{2^i+1}tr(x/a)tr((x/a)^{2^i+1})=x^{2^i+1}+(ax^{2^i}+a^{2^i}x
+a^{2^i+1})(tr(x/a)+tr((x/a)^{2^i+1})).$$
Therefore,
{\setlength\arraycolsep{2pt}
\begin{eqnarray*}
F_2\circ F_1^{-1}(x)&=&(x+a\ tr(x/a)+a\ 
tr((x/a)^{2^i+1}))^{2^i+1}+a^{2^i+1}tr(a^{-(2^i+1)}(x+a\ tr(x/a)\\
&+&a\ tr((x/a)^{2^i+1}))^{2^i+1})+a^{2^i+1}tr(a^{-1}(x+a\ tr(x/a)+a\ 
tr((x/a)^{2^i+1})))\\
&=&x^{2^i+1}+(ax^{2^i}+a^{2^i}x+a^{2^i+1})(tr(x/a)+tr((x/a)^{2^i+1}))+a^{2^i+1}tr((x/a)^{2^i}\\
&+&(x/a)+1)(tr(x/a)+tr((x/a)^{2^i+1}))+2a^{2^i+1}tr(x/a)+2a^{2^i+1}tr((x/a)^{2^i+1})\\
&=&a^{2^i+1}[(x/a)^{2^i+1}+((x/a)^{2^i}+(x/a))tr((x/a)+(x/a)^{2^i+1})]=a^{2^i+1}F^\prime(x/a),
\end{eqnarray*}}
where $F^{\prime}(x)=x^{2^i+1}+(x^{2^i}+x)tr(x+x^{2^i+1})$.
By Proposition \ref{prop1} the function $F^{\prime}$ is AB.\\ 
 It is not difficult to see that the algebraic degree of $F'$ is 3 for $m>3$.
 On the other hand
$tr(F^\prime(x))=tr(x^{2^i+1})$ and
$d^{\circ}(tr(F^\prime(x)))=2$. It follows from Proposition \ref{prop5} that 
$F^\prime$
is EA-inequivalent to any power function.\hfill $\Box$\\

\noindent \textbf{Remark 4 } It was conjectured in \cite{CC-carchzin} that any 
AB function is  EA-equivalent to a permutation. The AB function from Theorem 1 
serves as a counterexample for this conjecture.
It was checked by the help of a computer, that for no linear function $L$ on
$\mathbb{F}_{2^5}$ the sum $F^\prime+L$ is  a permutation for the AB
function $F^\prime(x)=x^{2^i+1}+(x^{2^i}+x)tr(x^{2^i+1}+x)$, $\gcd(5,i)=1$. 
Thus, $F^\prime$ is  EA-inequivalent to any permutation but it is CCZ-equivalent 
to the permutation $x^{2^i+1}$.
For $m$ even the existence of APN permutations on $\mathbb{F}_{2^m}$ is an 
open problem. It is shown by Nyberg that there exist no quadratic APN permutations when $m$ is even \cite{nyb2}. But as we have seen the nonexistence of permutations EA-equivalent to the Gold functions does not mean yet that there exist no permutations CCZ-equivalent to them.
\hfill $\Diamond$
\begin{Theorem}\label{thm2}
The function $F^\prime:\mathbb{F}_{2^m}\to \mathbb{F}_{2^m}$, $m\ge 4$ even,
\begin{displaymath}
F^\prime(x)=x^{2^i+1}+(x^{2^i}+x+1)tr(x^{2^i+1}), \textbf{  } \gcd(m,i)=1,
\end{displaymath}
is  an APN function which is EA-inequivalent to any power function.
\end{Theorem}
\emph{Proof.} For $a \in \mathbb{F}_{2^m}^{*},$ the linear mapping $${\cal 
L}(x,y)=(L_1,L_2)(x,y)=(x+a\ tr(a^{-(2^i+1)}y),y)$$ is an involution and, 
obviously, ${\cal L}(V)=V^\prime$. Thus, by Proposition \ref{prop2} the function 
$F$ is CCZ-equivalent to $F_2\circ F_1^{-1}$, where
{\setlength\arraycolsep{2pt}
\begin{eqnarray*}
F_1(x)&=&L_1(x,F(x))= x+a\ tr((x/a)^{2^i+1}),\\
F_2(x)&=&L_2(x,F(x))= x^{2^i+1}.
\end{eqnarray*}}
The function $F_1$ is an involution:
{\setlength\arraycolsep{2pt}
\begin{eqnarray*}
F_1\circ F_1(x)&=&x+a\ tr((x/a)^{2^i+1})+a\ 
tr(a^{-(2^i+1)}(x^{2^i+1}+ax^{2^i}tr((x/a)^{2^i+1})\\
&+&a^{2^i}x\ tr((x/a)^{2^i+1})+a^{2^i+1}tr((x/a)^{2^i+1})))\\
&=&x+2a\ tr((x/a)^{2^i+1})+a\ tr((x/a)^{2^i}+(x/a)+1)tr((x/a)^{2^i+1})=x.
\end{eqnarray*}}
We have
{\setlength\arraycolsep{2pt}
\begin{eqnarray*}
F_2\circ F_1^{-1}(x)&=&(x+a\ 
tr((x/a)^{2^i+1}))^{2^i+1}=x^{2^i+1}+ax^{2^i}tr((x/a)^{2^i+1})\\
&+&a^{2^i}x\ tr((x/a)^{2^i+1})+a^{2^i+1}tr((x/a)^{2^i+1})\\
&=&a^{2^i+1}[(x/a)^{2^i+1}+((x/a)^{2^i}+(x/a)+1)tr((x/a)^{2^i+1})].
\end{eqnarray*}}
The function $F_2\circ F_1^{-1}$ is EA-equivalent to the function
$$F^\prime(x)=x^{2^i+1}+(x^{2^i}+x+1)tr(x^{2^i+1}).$$
Hence, $F^\prime$ is CCZ-equivalent to $F$ and it is APN by Proposition 
\ref{prop1}.
The algebraic degree of $F^\prime$ is $3$ and
$d^\circ(tr(F^\prime))=2$ as $tr(F^\prime(x))=tr(x^{2^i+1})$. Therefore,  
$F^\prime$ is  EA-inequivalent to power functions
by Proposition \ref{prop5}. \hfill $\Box$\\

Note that the proofs of Theorems \ref{thm1} and \ref{thm2} do not depend on the 
condition $\gcd (i,m)=1$. When $\gcd (i,m)=s$ then the functions $F'$ have the 
same differential and linear properties as $x^{2^i+1}$ and, therefore, if $m/s$ 
is odd they can be considered as the first polynomials with three valued Walsh 
spectrum $\{0,\pm 2^\frac{m+s}{2}\}$, which are EA-inequivalent to power 
functions. 

\begin{Theorem}\label{thm3}
The function $F^\prime:\mathbb{F}_{2^m}\to \mathbb{F}_{2^m}$, $m$  divisible by 
$6$,  
$$F^\prime(x)=[x+tr_{m/3}(x^{2(2^i+1)}+x^{4(2^i+1)})+tr(x)tr_{m/3}(x^{2^i+1}+x^{2^{2i}(2^i+1)})]^{2^i+1},$$
with $\gcd(m,i)=1$,  is  an APN function.
\end{Theorem}
\emph{Proof.} The linear function ${\cal L}:\mathbb{F}_{2^m}^2\to 
\mathbb{F}_{2^m}^2$, 
$${\cal L}(x,y)=(L_1,L_2)(x,y)=(x+tr_{m/3}(y^2+y^4),y)$$
is a permutation since it has the kernel $\{(0,0)\}$. For $m$ divisible by $6$, 
the function 
$$F_1(x)=L_1(x,F(x))=x+tr_{m/3}(x^{2(2^i+1)}+x^{4(2^i+1)})$$ 
is a permutation. Indeed, by Proposition \ref{prop3} the function $F_1$ is a 
permutation  
if for every $v\in \mathbb{F}_{2^m}$ such that $tr(v)=0$ and every $u\in  
\mathbb{F}_{2^m}^*$ the condition $tr_{m/3}((u^{2^i+1}v)^2+(u^{2^i+1}v)^4)\ne u$ 
holds. Obviously, for any $u\notin  \mathbb{F}_{2^3}$  the condition is 
satisfied and when $u\in  \mathbb{F}_{2^3}^*$ the condition 
$tr_{m/3}((u^{2^i+1}v)^2+(u^{2^i+1}v)^4)\ne u$ is equivalent to 
$(u^{2^i+1}tr_{m/3}(v))^2+(u^{2^i+1}tr_{m/3}(v))^4\ne u$.
Therefore, $F_1$ is a  permutation if,  for every $u,w\in \mathbb{F}_{2^3}^*$, 
$tr_3(w)=0$ the condition  $(u^{2^i+1}w)^2+(u^{2^i+1}w)^4\ne u$ is satisfied and 
that was easily checked by a computer. 

We show below that $F_{1}^{-1}=F_1\circ F_1\circ F_1\circ F_1\circ F_1$.\\
We denote 
$$T(x)=tr_{m/3}(x^{2^i+1}),$$
then
$$F_1(x)=x+T(x)^2+T(x)^4.$$
Since every element of $\mathbb{F}_8$ is equal to its $8$ power and the function 
$tr_{m/3}(x)$ is $0$ on $\mathbb{F}_8$, then
{\setlength\arraycolsep{2pt}
\begin{eqnarray*}
&&T\circ F_1(x)=tr_{m/3}[(x+tr_{m/3}(x^{2(2^i+1)}+x^{4(2^i+1)}))^{2^i+1}]\\
&=&tr_{m/3}[(x+tr_{m/3}(x^{2(2^i+1)}+x^{4(2^i+1)}))(x^{2^i}+tr_{m/3}(x^{2^{i+1}(2^i+1)}\\
&+&x^{2^{i+2}(2^i+1)}))]=tr_{m/3}[x^{2^i+1}+x\ 
tr_{m/3}(x^{2^{i+1}(2^i+1)}+x^{2^{i+2}(2^i+1)})\\
&+&x^{2^i}tr_{m/3}(x^{2(2^i+1)}+x^{4(2^i+1)})]=tr_{m/3}(x)tr_{m/3}(x^{2^{s+1}(2^i+1)}+x^{2^{s+2}(2^i+1)})\\
&+&tr_{m/3}(x^{2^s})tr_{m/3}(x^{2(2^i+1)}+x^{4(2^i+1)})+tr_{m/3}(x^{2^i+1}),
\end{eqnarray*}} 
where $s\equiv i\mod 3$.\\
Therefore,
{\setlength\arraycolsep{2pt}
\begin{eqnarray*}
&&F_1\circ F_1(x)=F_1(x)+[T\circ F_1(x))]^2+[T\circ F_1(x))]^4\\
&=&x+2tr_{m/3}(x^{2(2^i+1)}+x^{4(2^i+1)})+tr_{m/3}(x^2)tr_{m/3}(x^{2^{s+2}(2^i+1)}\\
&+&x^{2^{s}(2^i+1)})+tr_{m/3}(x^{2^{s+1}})tr_{m/3}(x^{4(2^i+1)}+x^{2^i+1})\\
&+&tr_{m/3}(x^4)tr_{m/3}(x^{2^{s}(2^i+1)}+x^{2^{s+1}(2^i+1)})+tr_{m/3}(x^{2^{s+2}})tr_{m/3}(x^{2^i+1}+x^{2(2^i+1)})
\end{eqnarray*}} 
Considering separately the cases $s=1$ and $s=2$ we get
{\setlength\arraycolsep{2pt}
\begin{eqnarray*}
F_1\circ F_1(x)&=&x+tr_{m/3}(x+x^2+x^4)tr_{m/3}(x^{2^i+1}+x^{2^s(2^i+1)})\\
&=&x+tr(x)tr_{m/3}(x^{2^i+1}+x^{2^s(2^i+1)})=x+(T(x)+T(x)^{2^s})tr(x).
\end{eqnarray*}} 
Like this we get
$$F_1\circ F_1\circ F_1\circ F_1\circ F_1 
(x)=x+T(x)^2+T(x)^4+tr(x)(T(x)+T(x)^{2^{2s}}),$$
$$F_1\circ F_1\circ F_1\circ F_1\circ F_1\circ F_1 (x)= x.$$
Thus,
{\setlength\arraycolsep{2pt}
\begin{eqnarray*}
F^\prime(x)&=&F_2 \circ 
F_1^{-1}(x)=[x+T(x)^2+T(x)^4+tr(x)(T(x)+T(x)^{2^{2s}})]^{2^i+1}=x^{2^i+1}+T(x)^{2^s+1}\\
&+&tr(x^{2^i+1})T(x)^{2^{2s}}+tr(x)(T(x)+T(x)^4)+xtr(x)(T(x)+T(x)^{2^s})\\
&+&x^{2^i}tr(x)(T(x)+T(x)^{2^{2s}})+x(T(x)+T(x)^{2^{2s}})+x^{2^i}(T(x)^2+T(x)^4),
\end{eqnarray*}} 
where $F_2(x)=L_2(x,F(x))=x^{2^i+1}.$

The function $F'$ has the algebraic degree 4, this can be shown by lengthy but 
 uncomplicated calculations. Hence, $F^\prime$ is  EA-inequivalent to other known 
 APN functions since for $m$ divisible by 6 we have no known APN functions of 
 algebraic degree 4.
\hfill  $\Box$\\

Let $F:\mathbb{F}_{2^m}\to \mathbb{F}_{2^m}$, $F(x)=x^{2^i+1}$ with 
$\gcd(m,i)=1$, $m$ odd  and $n$  a divisor of $m$. Then $A_F\cap 
V^\prime=\emptyset$ for the subgroup 
$$V^\prime=\{(b,x) : b\in\mathbb{F}_{2^n},  x\in\mathbb{F}_{2^m},  tr_{m/n}(x)=0 
\},$$
since if $(b,x)\in V^\prime$ then 
$tr(b^{-(2^i+1)}x)=tr_n(b^{-(2^i+1)}tr_{m/n}(x))=0$ and $(b,x)\notin A_F$.
Hence, $G_F$ is transversal to $V^\prime$. Using the subgroup $V^\prime$ we 
construct an AB function $F^\prime$ given in the following theorem.
\begin{Theorem}\label{thm5}
The function $F^\prime:\mathbb{F}_{2^m}\to \mathbb{F}_{2^m}$, where $m$  is odd 
and divisible by $n$, $m\ne n$ and $\gcd(m,i)=1$,
 {\setlength\arraycolsep{2pt}
\begin{eqnarray*} 
F^\prime(x)&=&x^{2^i+1}+tr_{m/n}(x^{2^i+1})+x^{2^i}tr_{m/n}(x)+x\ 
tr_{m/n}(x)^{2^i}\\
&+&[tr_{m/n}(x)^{2^i+1}+tr_{m/n}(x^{2^i+1})+tr_{m/n}(x)]^{\frac{1}{2^i+1}}(x^{2^i}+tr_{m/n}(x)^{2^i}+1)\\
&+&[tr_{m/n}(x)^{2^i+1}+tr_{m/n}(x^{2^i+1})+tr_{m/n}(x)]^{\frac{2^i}{2^i+1}}(x+tr_{m/n}(x)),
\end{eqnarray*}} 
is  an AB function which is EA-inequivalent to any power function.
\end{Theorem}
\emph{Proof.} Let $m$ be odd and divisible by $n$. Obviously, the linear 
function 
$${\cal L}(x,y)=(L_1,L_2)(x,y)=(x+tr_{m/n}(x)+tr_{m/n}(y),y+tr_{m/n}(x))$$
is a permutation on $\mathbb{F}_{2^m}^2$ and
$${\cal L}^{-1}(x,y)=(x+tr_{m/n}(y),y+tr_{m/n}(x)+tr_{m/n}(y)).$$
We have
 $$F_1(x)=x+tr_{m/n}(x)+tr_{m/n}(x^{2^i+1}),$$ 
  $$F_2(x)=x^{2^i+1}+tr_{m/n}(x). $$
By Proposition \ref{prop2}, the function $F_1$ is a permutation  since ${\cal 
L}^{-1}(V)=V'$ and $G_F$ is transversal to $V^\prime$.
We need the inverse of the function $F_1$ to construct $F^\prime=F_2\circ 
F_1^{-1}$.\\
For any fixed element $x\in \mathbb{F}_{2^m}$ we have
$$y=x+tr_{m/n}(x)+tr_{m/n}(x^{2^i+1})=x+u,$$
for some $u\in \mathbb{F}_{2^n}$, and, therefore, $x=y+u$. Then
$$y=(y+u)+tr_{m/n}(y+u)+tr_{m/n}((y+u)^{2^i+1})$$ which yields 
\begin{equation} \label{eq1}
u^{2^i+1}+u^{2^i}tr_{m/n}(y)+u(tr_{m/n}(y))^{2^i}+tr_{m/n}(y^{2^i+1})+tr_{m/n}(y)=0.
\end{equation} 
If $tr_{m/n}(y)\ne 0$  then we denote $v=u/tr_{m/n}(y)$ and  we get
$$v^{2^i+1}+v^{2^i}+v+\frac{tr_{m/n}(y^{2^i+1})+tr_{m/n}(y)}{(tr_{m/n}(y))^{2^i+1}}=0.$$
Since $v^{2^i+1}+v^{2^i}+v=(v+1)^{2^i+1}+1$ then 
$$v+1=\left[\frac{tr_{m/n}(y^{2^i+1})+tr_{m/n}(y)}{(tr_{m/n}(y))^{2^i+1}}+1\right]^{\frac{1}{2^i+1}}.$$
Replacing $v$ by $u/tr_{m/n}(y)$ we have
$$u=[(tr_{m/n}(y))^{2^i+1}+tr_{m/n}(y^{2^i+1})+tr_{m/n}(y)]^{\frac{1}{2^i+1}}+tr_{m/n}(y).$$
If $tr_{m/n}(y)=0$ then from the equality  (\ref{eq1}) we get 
$u=[tr_{m/n}(y^{2^i+1})]^{\frac{1}{2^i+1}}$ and we observe that $u$ equals again 
$[(tr_{m/n}(y))^{2^i+1}+tr_{m/n}(y^{2^i+1})+tr_{m/n}(y)]^{\frac{1}{2^i+1}}+tr_{m/n}(y).$
Thus, in all cases, we have 
$$F_1^{-1}(y)=y+u=y+[(tr_{m/n}(y))^{2^i+1}+tr_{m/n}(y^{2^i+1})+tr_{m/n}(y)]^{\frac{1}{2^i+1}}+tr_{m/n}(y)$$
and
 {\setlength\arraycolsep{2pt}
\begin{eqnarray*} 
F^\prime(x)&=&F_2\circ  
F_1^{-1}(x)=[x+[(tr_{m/n}(x))^{2^i+1}+tr_{m/n}(x^{2^i+1})+tr_{m/n}(x)]^{\frac{1}{2^i+1}}+tr_{m/n}(x)]^{2^i+1}\\
&+&tr_{m/n}[x+[(tr_{m/n}(x))^{2^i+1}+tr_{m/n}(x^{2^i+1})+tr_{m/n}(x)]^{\frac{1}{2^i+1}}+tr_{m/n}(x)]\\
&=& 
x^{2^i+1}+tr_{m/n}(x^{2^i+1})+tr_{m/n}(x)+x^{2^i}tr_{m/n}(x)+x(tr_{m/n}(x))^{2^i}\\
&+&[(tr_{m/n}(x))^{2^i+1}+tr_{m/n}(x^{2^i+1})+tr_{m/n}(x)]^{\frac{1}{2^i+1}}(x^{2^i}+(tr_{m/n}(x))^{2^i}+1)\\ 
&+&[(tr_{m/n}(x))^{2^i+1}+tr_{m/n}(x^{2^i+1})+tr_{m/n}(x)]^{\frac{2^i}{2^i+1}}(x+tr_{m/n}(x)).
 \end{eqnarray*}}
 
We show below that  
the function $F^\prime$ has the algebraic degree $n+2$. 
It means that the number of functions CCZ-equivalent to a Gold AB function and 
EA-inequivalent to it is not smaller than the number of divisors of $m$.

The inverse of $x^{2^i+1}$ on $F_{2^n}$ is $x^d$, where 
$$d=\sum_{k=0}^{\frac{n-1}{2}}2^{2ik},$$ 
and $x^d$ has the algebraic degree $\frac{n+1}{2}$ \cite{nyb}.
Obviously, 
$((tr_{m/n}(x))^{2^i+1}+tr_{m/n}(x^{2^i+1})+tr_{m/n}(x))^d$ 
has the algebraic degree $n+1$ if and only if 
$((tr_{m/n}(x))^{2^i+1}+tr_{m/n}(x^{2^i+1}))^d$ has this algebraic degree.

We assume that $m\ne n$ and $n\ne 1$ since when $m=n$ we get  
$F^\prime(x)=x^{\frac{1}{2^i+1}}+x$ and Theorem \ref{thm1} gives the case $n=1$.

We have
$$tr_{m/n}(x)=\sum_{k=0}^{\frac{m}{n}-1}x^{2^{kn}}$$
and
$$(tr_{m/n}(x))^{2^i+1}+tr_{m/n}(x^{2^i+1})=\sum_{k=0}^{\frac{m}{n}-1}x^{2^{kn}}\sum_{k=0}^{\frac{m}{n}-1}x^{2^{kn+i}}+\sum_{k=0}^{\frac{m}{n}-1}(x^{2^i+1})^{2^{kn}}$$
$$=\sum_{k,j=0}^{\frac{m}{n}-1}x^{2^{kn}+2^{jn+i}}+\sum_{k=0}^{\frac{m}{n}-1}x^{2^{kn}+2^{kn+i}}=\sum_{\begin{subarray}{c} 
k,j=0\\k\ne j\end{subarray}}^{\frac{m}{n}-1}x^{2^{kn}+2^{jn+i}}.$$

Note that we have
$$[(tr_{m/n}(x))^{2^i+1}+tr_{m/n}(x^{2^i+1})]^{2^{2i}+1}=\sum_{\begin{subarray}{c} 
k,j=0\\k\ne 
j\end{subarray}}^{\frac{m}{n}-1}x^{2^{kn}+2^{jn+i}}\sum_{\begin{subarray}{c} 
k,j=0\\k\ne j\end{subarray}}^{\frac{m}{n}-1}x^{2^{kn+2i}+2^{jn+3i}}$$
$$=\sum_{\begin{subarray}{c} k,j,s,t=0\\k\ne j, s\ne t 
\end{subarray}}^{\frac{m}{n}-1}x^{2^{kn}+2^{jn+i}+2^{sn+2i}+2^{tn+3i}}.$$
Similarly, we have
\begin{equation}\label{1}
((tr_{m/n}(x))^{2^i+1}+tr_{m/n}(x^{2^i+1}))^d=\sum_{(k_0,...,k_n)\in 
I}x^{\sum_{s=0}^{n}2^{k_sn+si}}, 
\end{equation}
where $I=\{(l_0,...,l_n):\quad 0\le l_t\le \frac{m}{n}-1,\quad l_{2t}\ne  
l_{2t+1}\}$.

The equality $k_sn+si= k_tn+ti$ is possible for $0\le s<t\le n$ only for $s=0$ 
and $t=n$. Indeed, if $k_sn+si= k_tn+ti$ then $(k_s-k_t)n=(t-s)i$. Since 
$\gcd(n,i)=1$ and $0\le t,s\le n$ then $t=n,$ $s=0$ and $k_s-k_t=i$.

For simplicity we consider now the equality (\ref{1}) in the case $i=1$. In the sum 
$\sum_{s=0}^{n}2^{k_sn+s}$ the largest possible item is 
$2^{(\frac{m}{n}-1)n+n}=2^m$. Therefore, when $k_n \ne \frac{m}{n}-1$  the sum 
is smaller than $2^m-1$. Besides, all items in the sum are different modulo 
$2^m-1$ except the case when $k_0=0$ and $k_n = \frac{m}{n}-1$ and in the cases 
where $k_0= k_n+1$.
Therefore, when  $k_0=k_n=1$ the number 
 $\sum_{s=0}^{n}2^{k_sn+s}$ has the weight $n+1$. On the other hand, when 
 $k_0=k_n=1$ and $k_1=k_{n-1}=0$ 
the term $$x^{\sum_{s=0}^{n}2^{k_sn+s}}$$ does not vanish in (\ref{1}). Indeed, 
if 
$${\sum_{s=0}^{n}2^{k_sn+s}}\equiv {\sum_{p=0}^{n}2^{t_pn+p}} \mod (2^m-1)$$
then we have only two possibilities:\\
1) for any $s$ there exists $p$ such that $k_sn+s=t_pn+p$ (and vice versa). Then 
$(k_s-t_p)n=p-s$ and since $0\le s,p\le n$ then $k_0=t_n+1$, 
$t_0=k_n+1$ and $k_s=t_s$ for $s\ne 0,n$.  If $k_0=k_n=1$ then $t_n=0$, $t_0=2$. 
But in our case $k_0=k_n=1$, $t_1=k_1=0$, $t_{n-1}=k_{n-1}=0$ and, therefore, 
$t_n\ne 0$ since $t_{n-1}\ne t_n$.
\\
2) if $t_n=\frac{m}{n}-1$ then $k_0$ must be equal to $0$ or 
$k_n=\frac{m}{n}-1$, but $k_0=k_n=1$.

Thus, when $k_0=k_n=1$ and $k_1=k_{n-1}=0$ (for permissible $k_s$, $1<s<n-1$) 
the term $$x^{\sum_{s=0}^{n}2^{k_sn+s}}$$ has the algebraic degree $n+1$ and it 
does not vanish in (\ref{1}). 

If $n\ge 5$ we can also take 
 $k_2=1,k_3=k_4=0$ and then we get 
 \begin{equation}\label{2} 
 \sum_{s=0}^{n}2^{k_sn+s}=2^n+2+2^{n+2}+2^3+2^4+...+2^{2n}. 
 \end{equation} 
 
We have 
 $$((tr_{m/n}(x))^3+tr_{m/n}(x^3))^d(x^2+tr_{m/n}(x)^{2})+((tr_{m/n}(x))^3+tr_{m/n}(x^3))^{2d}(x+tr_{m/n}(x))$$ 
 $$ 
 =\sum_{(k_0,...,k_n)\in I}x^{\sum_{s=0}^{n}2^{k_sn+s}}\sum_{1\le k\le 
 m/n-1}x^{2^{nk+1}}+\sum_{(k_0,...,k_n)\in 
 I}x^{\sum_{s=0}^{n}2^{k_sn+s+1}}\sum_{1\le j\le m/n-1}x^{2^{nj}} 
 $$ 
\begin{equation}\label{*}
=\sum_{\begin{subarray}{c}(k_0,...,k_n)\in I\\ 1 \le k\le 
 m/n-1\end{subarray}}x^{2^{nk+1}+\sum_{s=0}^{n}2^{k_sn+s}}+\sum_{\begin{subarray}{c}(k_0,...,k_n)\in I\\ 1 \le j\le m/n-1\end{subarray}}x^{2^{nj}+\sum_{s=0}^{n}2^{k_sn+s+1}}.
\end{equation}  
 We consider the item with the exponent 
 \begin{equation}\label{3} 
 2^n+2+2^{n+2}+2^3+2^4+...+2^{2n}+2^{nk+1} 
 \end{equation} 
 from the first sum in (\ref{*}). It is easy to see that 
 $2^n+2+2^{n+2}+2^3+2^4+...+2^{2n}+2^{nk+1}<2^m$ since $k\le m/n-1$. In this sum 
 $nk+1=k_sn+s$ only if $s=1$. But then $k=k_1=0$ which is in contradiction with 
 $1\le k$. Thus, the number given by this sum has the weight $n+2$. The item with 
 the exponent (\ref{3}) does not vanish. Indeed, if there is another item in the 
 first sum of (\ref{*}) with this exponent then 
 \begin{equation*} 
 2^n+2+2^{n+2}+2^3+2^4+...+2^{2n}+2^{nk+1}=2^{nj+1}+\sum_{s=0}^{n}2^{k_sn+s}. 
 \end{equation*} 
 If $k= j$ then (\ref{2}) is equal to another sum $\sum_{s=0}^{n}2^{k_sn+s}$ and 
 we already showed that it is impossible. If $k\ne j$ then $k_1=k$ and $j=0$ 
 while $1\le j$.\\ 
 Assume there exists an item in the second sum of (\ref{*}) with the exponent (\ref{3}) 
 then 
 \begin{equation*} 
 2^n+2+2^{n+2}+2^3+2^4+...+2^{2n}+2^{nk+1}=2^{nj}+\sum_{s=0}^{n}2^{k_sn+s+1} 
 \end{equation*} 
 for some $j,$  $1\le j\le m/n-1$ and $(k_0,...,k_n)\in I$. We have $3=k_sn+s+1 
 \mod m$ for some $s$, $0\le s\le n$. Then $k_sn=2-s$ or $k_sn=m-(s-2)$ and this 
 is possible only if $k_2=0$ or $k_2=m/n$, but since $0\le k_s\le m/n-1$, then 
 $3=k_sn+s+1 \mod m$ only if $k_2=0$. The same arguments show that $4=k_sn+s+1 
 \mod m$ only if $k_3=0$ and that is in contradiction with the condition 
 $k_{2t}\ne k_{2t+1}$. Therefore the item  with the exponent (\ref{3}) does not 
 vanish in (\ref{*}) and then it does not vanish in the sum presenting the function $F'$.  This completes the proof that $F'$ has the algebraic degree $n+2$.
   
The algebraic degree of the function $tr(F^\prime(x))$ is not greater than $n+1$ since
{\setlength\arraycolsep{2pt}
\begin{eqnarray*} 
tr(F^\prime(x))&=&tr(x^{2^i+1}+tr_{m/n}(x^{2^i+1})+tr_{m/n}(x)+x^{2^i}tr_{m/n}(x)+x(tr_{m/n}(x))^{2^i})\\
&+&tr_n([(tr_{m/n}(x))^{2^i+1}+tr_{m/n}(x^{2^i+1})+tr_{m/n}(x)]^{\frac{1}{2^i+1}}tr_{m/n}(x^{2^i}+(tr_{m/n}(x))^{2^i}+1))\\
&+&tr_n([(tr_{m/n}(x))^{2^i+1}+tr_{m/n}(x^{2^i+1})+tr_{m/n}(x)]^{\frac{2^i}{2^i+1}}tr_{m/n}(x+tr_{m/n}(x)))\\
&=&tr(x)+tr(x^{2^i}tr_{m/n}(x))+tr(x(tr_{m/n}(x))^{2^i})\\
&+&tr([(tr_{m/n}(x))^{2^i+1}+tr_{m/n}(x^{2^i+1})+tr_{m/n}(x)]^{\frac{1}{2^i+1}}).
\end{eqnarray*}}
On the other hand $d^\circ (tr(F'(x)))$ is not 0 or 1. Indeed, for any  $x\in \mathbb{F}_{2^n}$ we have
$$tr(F^\prime(x))=tr(x)+2tr(x^{2^i+1})+tr([2x^{2^i+1}+x]^{\frac{1}{2^i+1}})=tr(x)+tr(x^{\frac{1}{2^i+1}}).$$
The function $tr(x^{\frac{1}{2^i+1}})$ has the algebraic degree $(n+1)/2$. Indeed, $d^\circ (tr(x^{\frac{1}{2^i+1}}))\in\{0,d^\circ (x^{\frac{1}{2^i+1}})\}$ and since $x^{\frac{1}{2^i+1}}$ is a permutation then $tr(x^{\frac{1}{2^i+1}})$ is not a constant and $d^\circ (tr(x^{\frac{1}{2^i+1}}))=d^\circ (x^{\frac{1}{2^i+1}})=(n+1)/2$.
Hence, $tr(F'(x)$ is not linear on $\mathbb{F}_{2^n}$ and then it cannot be linear on $\mathbb{F}_{2^m}$.

Thus, the function $F^\prime$ is  EA-inequivalent to any power 
function by Proposition \ref{prop5}.\hfill  $\Box$\\

\section{Conclusion}
We have shown that there exist APN and AB functions which are EA-inequivalent to 
power functions, and therefore, which are new, up to EA-equivalence. But these 
functions are CCZ-equivalent to the Gold functions. We leave two open 
problems:\\
- finding classes of functions CCZ-equivalent to other known APN or AB 
functions, which would be  EA-inequivalent to all known APN and AB functions (or 
even, inequivalent to power functions);\\
- finding classes of APN or AB functions which would be CCZ-inequivalent to all 
known APN and AB functions (or even, CCZ-inequivalent to power functions).


\begin{thebibliography}{99}

\bibitem{BFF} T. Bending, D. Fon-Der-Flaass. Crooked functions, bent functions 
and distance-regular graphs. {\em Electron. J. Comb.}, 5(R34), 14, 1998.

\bibitem{B-D} T.~Beth and  C.~Ding. \newblock  On almost perfect nonlinear 
permutations. \newblock {\em  Advances in Cryptology-EUROCRYPT'93, Lecture Notes 
in Computer Science}, 765, Springer-Verlag, New York, pp. 65-76, 1993.

\bibitem{B-Sh} E.~Biham and A.~Shamir. \newblock Differential Cryptanalysis of 
{DES}-like Cryptosystems. \newblock {\em Journal of Cryptology}, vol. 4, No.1, 
pp. 3-72, 1991.

\bibitem{BCP} L.~Budaghyan, C.~Carlet, A.~Pott. New Constructions of Almost Perfect Nonlinear and Almost Bent Functions. \newblock {\em Proceedings of the Workshop on Coding and Cryptography 2005}, P.~Charpin and \O.~Ytrehus eds, pp. 306-315, 2005.

\bibitem{CCD3} A.~Canteaut, P.~Charpin and H.~Dobbertin. \newblock A new 
characterization of almost bent functions. \newblock {\em Fast Software 
Encryption 99, Lecture Notes in Computer Science} 1636, L. Knudsen edt, pp. 
186-200. Springer-Verlag, 1999.

\bibitem{CCD2} A.~Canteaut, P.~Charpin and H.~Dobbertin. \newblock  Binary 
$m$-sequences with three-valued crosscorrelation: A proof of Welch's conjecture. 
\newblock {\em IEEE Trans. Inform. Theory}, 46 (1), pp. 4-8, 2000.

\bibitem{CCD} A.~Canteaut, P.~Charpin, H.~Dobbertin. Weight divisibility of 
cyclic codes , highly nonlinear functions on $\mathbb{F}_{2^m}$, and 
crosscorrelation of maximum-length sequences. {\em SIAM Journal on Discrete 
Mathematics}, 13(1), pp. 105-138, 2000.

\bibitem{Cthesis} C.~Carlet. Codes de Reed-Muller,  codes de Kerdock et de 
Preparata, PhD thesis, Publication of LITP, Institut Blaise Pascal, 
Universit{\'e} Paris 6, 90.59, 1990.

\bibitem{Cbook1} C.~Carlet. \newblock Boolean Functions for Cryptography and 
Error Correcting Codes. \newblock Chapter of the monography {\em Boolean Methods 
and 
 Models}, Y.~Crama and 
 P.~Hammer eds,  Cambridge University Press, to appear (winter 2005-2006).

\bibitem{Cbook} C.~Carlet. \newblock  Vectorial Boolean Functions for 
Cryptography. \newblock Chapter of the monography {\em Boolean Methods and 
 Models}, Y.~Crama and 
 P.~Hammer eds,  Cambridge University Press, to appear (winter 2005-2006).

\bibitem{CC-carchzin}
C.~Carlet, P.~Charpin and V.~Zinoviev.
\newblock Codes, bent functions and permutations suitable for {DES}-like
  cryptosystems.
\newblock {\em Designs, Codes and Cryptography}, 15(2), pp. 125-156, 1998.

\bibitem{CD}  C.~Carlet and C.~Ding.  \newblock Highly Nonlinear Mappings. 
\newblock {\em Special Issue "Complexity Issues in Coding and Cryptography", 
dedicated to Prof. Harald Niederreiter on the occasion of his 60th birthday, 
Journal of Complexity} 20, pp. 205-244, 2004.

\bibitem{CV} F. Chabaud and S. Vaudenay. Links between differential and linear 
cryptanalysis, {\em Advances in Cryptology -EUROCRYPT'94, Lecture Notes in 
Computer Science}, Springer-Verlag, New York, 950, pp. 356-365, 1995.

\bibitem{CourtPier} N. Courtois and J. Pieprzyk.
Cryptanalysis of block ciphers with overdefined systems of equations.
{\em Advances in cryptology--ASIACRYPT 2002,
Lecture Notes in Computer Science} 2501,
pp. 267-287, Springer, 2003.

\bibitem{D4} H.~Dobbertin. \newblock  One-to-One Highly Nonlinear Power 
Functions on $GF(2^n)$. \newblock {\em Appl. Algebra Eng. Commun. Comput.} 9 
(2), pp. 139-152, 1998.

\bibitem{D1} H.~Dobbertin. \newblock  Almost perfect nonlinear power functions 
over $GF(2^n)$: the Niho case. \newblock {\em Inform. and Comput.}, 151, pp. 
57-72, 1999.

\bibitem{D2} H.~Dobbertin. \newblock Almost perfect nonlinear power functions 
over $GF(2^n)$: the Welch case.  \newblock {\em IEEE Trans. Inform. Theory}, 45, 
pp. 1271-1275, 1999.

\bibitem{D3} H.~Dobbertin. \newblock Almost perfect nonlinear power functions 
over $GF(2^n)$: a new case for $n$ divisible by 5. \newblock D.~Jungnickel and 
H.~Niederreiter eds. {\em Proceedings of Finite Fields and Applications FQ5}, 
Augsburg, Germany, Springer, pp. 113-121, 2000.

\bibitem{Gold} R.~Gold. \newblock  Maximal recursive sequences with 3-valued 
recursive crosscorrelation functions. \newblock {\em IEEE Trans. Inform. 
Theory}, 14, pp. 154-156, 1968.


\bibitem{H-X} H.~Hollmann and Q.~Xiang. \newblock A proof of the Welch and Niho 
conjectures on crosscorrelations of binary $m$-sequences. \newblock {\em Finite 
Fields and Their Applications 7}, pp. 253-286, 2001.


\bibitem{J-W} H.~Janwa and R.~Wilson. \newblock Hyperplane sections of Fermat 
varieties in $P^3$ in char. 2 and some applications to cyclic codes. \newblock 
{\em Proceedings of AAECC-10, Lecture Notes in Computer Science}, vol. 673, 
Berlin, Springer-Verlag, pp. 180-194, 1993.

\bibitem{Kasami} T.~Kasami. \newblock  The  weight enumerators for several 
classes of subcodes of the second order binary Reed-Muller codes. \newblock {\em 
Inform. and Control}, 18, pp. 369-394, 1971.

\bibitem{LW} G.~Lachaud and J.~Wolfmann. \newblock The Weights of the 
Orthogonals of the Extended Quadratic Binary Goppa Codes. {\em IEEE Trans. 
Inform. Theory}, vol. 36, pp. 686-692, 1990.


\bibitem{M} M.~Matsui. \newblock Linear cryptanalysis  method for {DES} cipher. 
\newblock {\em Advances in Cryptology-EUROCRYPT'93, Lecture Notes in Computer 
Science}, Springer-Verlag, pp. 386-397, 1994.

\bibitem{nyb1} K.~Nyberg. \newblock On the construction of highly nonlinear 
permutations.
{\em Advances in Cryptography, EUROCRYPT'92, Lecture Notes in Computer Science}, 
Springer-Verlag, 658, pp. 92-98, 1993.

\bibitem{nyb} K.~Nyberg. \newblock Differentially uniform mappings for 
cryptography, {\em Advances in Cryptography, EUROCRYPT'93, Lecture Notes in 
Computer Science}, Springer-Verlag, New York, 765,  pp. 55-64, 1994.

\bibitem{nyb2} K.~Nyberg. \newblock S-boxes and Round Functions with 
Controllable Linearity and Differential Uniformity. \newblock  {\em Proceedings 
of Fast Software Encryption 1994, Lecture Notes in Computer Science} 1008, pp. 
111-130, 1995.

\bibitem{P} A.~Pott. \newblock Nonlinear functions in Abelian groups and 
relative difference sets. \newblock  {\em Discrete Applied Math.} 138, pp. 
177-193, 2004.

\bibitem{sid} V.~Sidelnikov. \newblock On mutual correlation of sequences, {\em 
Soviet Math. Dokl.}, 12(1971), pp. 197-201.

\bibitem{T} H. M. Trachtenberg. On the cross-correlation functions of maximal 
linear recurring sequences. PhD Thesis, University of Southern California, 1970.
\end{thebibliography}
\end{document}